\documentclass[reqno]{amsart}
\usepackage{hyperref}
\newtheorem*{theoA}{Theorem A}
\newtheorem*{theoB}{Theorem B}
\newtheorem*{theoC}{Theorem C}
\newtheorem*{theoD}{Theorem D}
\newtheorem*{theoE}{Theorem E}

\newtheorem{theo}{Theorem}[section]
\newtheorem{lem}{Lemma}[section]

\newtheorem{ques}{Question}[section]
\newtheorem{exm}{Example}[section]

\newtheorem{rem}{Remark}[section]
\newcommand{\ol}{\overline}
\newcommand{\be}{\begin{equation}}
\newcommand{\ee}{\end{equation}}
\newcommand{\beas}{\begin{eqnarray*}}
	\newcommand{\eeas}{\end{eqnarray*}}
\newcommand{\bea}{\begin{eqnarray}}
\newcommand{\eea}{\end{eqnarray}}

\numberwithin{equation}{section}

\begin{document}
\title[\hfilneg \hfil  Fermat-type partial differential-difference equations in $ \mathbb{C}^n $]
{Solutions of Fermat-type partial differential-difference equations in $ \mathbb{C}^n $}
\author[G. Haldar \hfil
\hfilneg]
{Goutam Haldar }

\address{Goutam Haldar, Department of Mathematics, Malda College, Malda - 732101, West Bengal, India.}
\email{goutamiit1986@gmail.com, goutamiitm@gmail.com}

%\thanks{Submitted   Aug.    16,  2017.}

\maketitle
%\noindent Uniqueness of difference-differential polynomials of meromorphic functions sharing a small function IM
%\setcounter{page}{288}

\begin{abstract}
For two meromorphic functions $ f $ and $ g $, the equation $ f^m+g^m=1 $ can be regarded as Fermat-type equations. Using Nevanlinna theory for meromorphic functions in several complex variables, the main purpose of this paper is to investigate the properties of the transcendental entire solutions of Fermat-type difference and partial differential-difference equations in $ \mathbb{C}^n $. In addition, we find the precise form of the transcendental entire solutions in $ \mathbb{C}^2 $ with finite order of the Fermat-type  partial differential-difference equation \beas \left(\frac{\partial f(z_1,z_2)}{\partial z_1}\right)^2+(f(z_1+c_1,z_2+c_2)-f(z_1,z_2))^2=1\eeas and \beas f^2(z_1,z_2)+P^2(z_1,z_2)\left(\frac{\partial f(z_1+c_1,z_2+c_2)}{\partial z_1}-\frac{\partial f(z_1,z_2)}{\partial z_1}\right)^2=1,\eeas where $P(z_1,z_2)$ is a polynomial in $\mathbb{C}^2$. Moreover, one of the main results of the paper significantly improved the result of Xu and Cao [Mediterr. J. Math. (2018) 15:227 , 1-14 and Mediterr. J. Math. (2020) 17:8, 1-4].
\end{abstract}

\section{\textbf{Introduction}}
This paper mainly deals with the existence of transcendental entire solutions with finite order of Fermat-type difference and partial differential-difference equations in several complex variables.  We adopt the standard notations of the Nevanlinna theory of meromorphic functions in one variable (see \cite{Hayman & 1964, Laine & 1993, Yi & Yang & 1995}). In the past several decades, considerable attention had been paid to the existence of entire solutions for Fermat-type equation $ x^n+y^m=1 $. Specially, in $1995$, Wiles and Taylor \cite{Tailor & Wiles & 1995, Wiles & Ann. Math. 1995} pointed out that the Fermat-type equation $x^n+y^m=1$, where $m,n\in \mathbb{N}$ does not admit nontrivial solution in rational numbers for $m=n\geq3$, and does exist nontrivial solution in rational numbers for $m=n=2$. Initially, Fermat-type functional equations were investigated by Gross \cite{Gross & Bull. Amer. & 1966, Gross & Amer. Math. & 1966}, Montel \cite{Montel & Paris & 1927}.\vspace{1.2mm}
\par In $1970$, Yang \cite{Yang & 1970} considered the following functional equation \bea\label{e1.1} f^n+g^m=1\eea and proved the following interesting result.
\begin{theoA}\cite{Yang & 1970}
There are no non-constant entire solutions of the functional equation $(\ref{e1.1})$ if $m$, $n$ are positive integers satisfying $1/m+1/n<1$.
\end{theoA}
\par In recent years, looking for non-constant entire as well as  meromorphic solutions of Fermat-type difference as well as differential-difference equations have been studied extensively after the development of the difference analogous lemma of the logarithmic derivative by Halburd and Korhonen \cite{Halburd & Korhonen & 2006, Halburd & Korhonen & Ann. Acad. & 2006}, and Chiang and Feng \cite{Chiang & Feng & 2008}, independently. As a result, successively a lot of investigations done by many researchers in this direction (see \cite{Chen & Gao & CKMS & 2015, Hal & Kor & lon & 2007, Hu & Yang & 2014, Li & Arch. Math. & 2008, Liu & JMAA & 2009, Liu & Cao & EJDE & 2013, Liu & Cao & Arch. Math. & 2012, Liu & Yang & 2016, Lu & Kodai & 2003, Tang & Liao & 2007, Wang & Xu & Tu & 2019, Yang & Li & 2004, Zhang & Liao & 2013}).\vspace{1.2mm}
\par In view of Theorem A, Liu \textit{et. al.} \cite{Liu & Cao & Arch. Math. & 2012} have showed that there is no finite order transcendental entire solution of Fermat-type difference equation $f^n(z)+f^m(z+c)=1$ when $n>m>1$ or $n=m>2$, and for the case $n=m=2$, $f(z)$ must be of the form  $f(z)=\sin(Az + B)$, where $c(\neq0)$, $B\in\mathbf{C}$ and $A=(4k+1)\pi/2c$, $k$ is an integer. Later, in $ 2013 $, Liu and Yang \cite{Liu & Yang & CMFT & 2013} extended this result by considering the Fermat-type difference equation $f^2(z)+P^2(z)f^2(z+c)=Q(z)$ where $P(z)$ and $Q(z)$ are two non-zero polynomials.\vspace{1.2mm}
\par In 2018,  Xu and Cao \cite{Xu & Cao & 2018} extended the above result to several complex variables as follows. 
\begin{theoB}\cite{Xu & Cao & 2018}
Let $c=(c_1,c_2,\ldots,c_n)\in \mathbb{C}^n\setminus\{0\}$. Then any non-constant entire solution $f:\mathbb{C}^n\rightarrow\mathbb{P}^{1}(\mathbb{C})$ with finite order of the Fermat-type difference equation $f(z)^2+f(z+c)^2=1$ has the form
of $f(z)=\cos(L(z)+B)$, where $L$ is a linear function of the form $L(z)=a_1z_1+\cdots+ a_nz_n$ on $\mathbb{C}^n$ such that $L(c)=-\pi/2-2k\pi$ $(k\in\mathbb{Z})$, and $B$ is a constant on $\mathbb{C}$.
\end{theoB}
It would be interesting of finding the solutions as well as the precise form of the solutions of the equation considered in Theorem B. Henceforth, motivated by the results of Xu and Cao \cite{Xu & Cao & 2018}, and Liu and Yang \cite{Liu & Yang & CMFT & 2013}, in this paper, we are mainly interested to generalize Theorem B by considering the following equation \bea\label{e2.1} f^2(z)+P^2(z)f^2(z+c)=Q(z),\eea where $P(z)$ and $Q(z)$ (with $Q(0)\neq 0$) are two non-zero polynomials in $\mathbb{C}^n$. It is easy to see that \eqref{e2.1} is a general setting of the equation $ f(z)^2+f(z+c)^2=1 $. Consequently, we prove the following result.
\begin{theo}\label{t1}
Let $f(z)$ be a transcendental entire function defined on $\mathbb{C}^n$ and $c=(c_1,c_2,\ldots,c_n)\in\mathbb{C}^n$. If $f(z)$ is a solution of \eqref{e2.1} such that $T(r,P)=o(T(r,f))$ and $T(r,Q)=o(T(r,f))$, then $P(z)=\pm1$, $Q(z)$ reduces to a constant and $f(z)$ is of the form $f(z)=\cos(L(z)+B)$ , where $L(z)=a_1z_1+\cdots+ a_nz_n$ such that $L(c)=-\pi/2-2k\pi$ $(k\in\mathbb{Z})$, and $B$ is a constant on $\mathbb{C}$.	
\end{theo}
\par Besides finding the solutions of Fermat-type difference equations, Fermat-type differential-difference equations are also studied by many researchers. For example, Xu and Cao \cite{Xu & Cao & 2018} have  investigated the entire solutions of Fermat-type partial differential-difference equation
\bea\label{e1.3} \left(\frac{\partial f(z_1,z_2)}{\partial z_1}\right)^n+f^m(z_1+c_1,z_2+c_2)=1\eea and obtained the following interesting result for the functions in $ \mathbb{C}^2 $.

\begin{theoC}\cite{Xu & Cao & 2018}
Let $c=(c_1,c_2)$ be a constant in $\mathbb{C}^2$. Then the Fermat-type partial differential-difference equation $(\ref{e1.3})$ does not have any transcendental entire solution with finite order, where $m$ and $n$ are two distinct positive integers.
\end{theoC}
In $2020$, Xu and Wang \cite{Xu & wang & 2020} generalized Theorem C by considering the following Fermat-type partial differential-difference equation \bea\label{e1.4} \left(\frac{\partial f(z_1,z_2)}{\partial z_1}+\frac{\partial f(z_1,z_2)}{\partial z_2}\right)^n+f^m(z_1+c_1,z_2+c_2)=1\eea
and proved the following result.
\begin{theoD}\cite{Xu & wang & 2020}
Let $c=(c_1,c_2)$ be a constant in $\mathbb{C}^2$ and $ m $, $ n $ be two positive integers. If the Fermat-type partial differential-difference equation \eqref{e1.4} satisfies one of the following conditions:
\begin{enumerate}
	\item[(i)] $ m>n $;
	\item[(ii)] $ n>m\geq 2 $, 
\end{enumerate}
then \eqref{e1.4} does not have any finite order transcendental entire solutions.
\end{theoD} 
As far as we know, there are few results about the complex difference and complex difference equations for the functions in several complex variables. Furthermore, it appears that the Fermat-type mixed partial differential-difference equations in several complex variables has not been addressed in the literature before. In order to generalize and also to establish a result which combines Theorem C and Theorem D, in our investigation, we consider the following partial differential-difference equation
\bea\label{e22.5} (\partial^{I}f(z_1,z_2)+\partial^{J}f(z_1,z_2))^n+f^m(z_1+c_1,z_2+c_2)=1,\eea where 
\beas \partial^{I}f(z_1,z_2)=\displaystyle\frac{\partial^{|I|}f(z_1,z_2)}{\partial z_1^{\alpha_1}\partial z_2^{\alpha_2}}\;\; \mbox{and}\;\; \partial^{J}f(z_1,z_2)=\displaystyle\frac{\partial^{|J|}f(z_1,z_2)}{\partial z_1^{\beta_1}\partial z_2^{\beta_2}} \eeas with $I=(\alpha_1,\alpha_2)$ and $J=(\beta_1,\beta_2)$ are two multi-index, where $\alpha_1,\alpha_2$, $\beta_1$ and $\beta_2$ are non-negative integers and $c_1, c_2\in\mathbb{C}$. We denote by $\mid I\mid$ to denote the length of $I$, that is, $\mid I\mid=\alpha_1+\alpha_2$. Similarly, for $J$ also. Throughout this paper, we denote $ z+w=(z_1+w_1, z_2+w_2) $ for any $ z=(z_1, z_2) $, $ w=(w_1, w_2) $ and $ c =(c_1, c_2) $.\vspace{1.2mm}
\par As a matter of fact, we prove the next result for any order Fermat-type partial differential-difference equation \eqref{e22.5}.
\begin{theo}\label{t2}
Let $c=(c_1,c_2)$ be a constant in $\mathbb{C}^2$ and $ m $, $ n $ be two positive integers. If the Fermat-type partial differential-difference equation \eqref{e22.5} satisfies one of the following conditions:
\begin{enumerate}
	\item[(i)] $ m>n $;
	\item[(ii)] $ n>m\geq 2 $, 
\end{enumerate}
then \eqref{e22.5} does not have any finite order transcendental entire solutions.
\end{theo}
\par In this regard, we recall here a result of Xu and Cao \cite{Xu & Cao & 2018} in which the authors have investigated the solutions of the following Fermat-type partial differential-difference equation
\bea\label{e2.3} \left(\frac{\partial f(z_1,z_2)}{\partial z_1}\right)^2+f^2(z_1+c_1,z_2+c_2)=1\eea and proved that the finite order transcendental entire solutions of \eqref{e2.3} must assume the form $f(z_1,z_2)=\sin(Az_1+B)$, where $A\in\mathbb{C}$ be such that $Ae^{iAc_1}=1$, and $B$ is a constant on $\mathbb{C}$ (see \cite[Theorem 1.2]{Xu & Cao & 2018}). The authors concluded that, all entire solutions $f(z_1,z_2)$ of \eqref{e2.3} must be a function of one variable $z_1$.
\begin{rem}
	 By the following example one can observe that the concluding part of the result of Xu and Cao is not true in general.
\end{rem}
\begin{exm}\label{exm1}
Let $f(z_1,z_2)=\sin(z_1+z_2)$ or $ \cos(z_1+z_2) $. Choose, $c=(c_1,c_2)\in \mathbb{C}^2$ such that $c_1+c_2=k\pi$ with $k\in \mathbb{Z}$. Then it is easy to see that $f(z_1,z_2)$ is a solution of the equation $(\ref{e2.3})$ but $ f $ is not a function of $ z_1 $ only.
\end{exm}
In 2020, Xu and Cao \cite{Xu & Cao & 2020} corrected \textbf{Theorem 1.2} of \cite{Xu & Cao & 2018} and obtained the following result.
\begin{theoE}\cite{Xu & Cao & 2020}
Let $c=(c_1,c_2)$ be a constant in $\mathbb{C}^2$. Then any transcendental entire solution with finite order of the Fermat-type partial differential-difference equation $(\ref{e2.3})$ has the form  $f(z_1,z_2)=\sin(Az_1+Bz_2+H(z_2))$, Where $A, B$ are constant on $\mathbb{C}$ satisfying $A^2=1$ and $Ae^{i(Ac_1+Bc_2)}=1$, and $H(z_2)$ is a polynomial in one variable $z_2$ such that $H(z_2)\equiv H(z_2 + c_2)$. In the special case whenever $c_2\neq0$, we have $f(z_1, z_2)=\sin (Az_1 + Bz_2 + Constant)$.
\end{theoE}
\par After a detailed study of the proof of result of Xu and Cao  \cite[Theorem 1.2]{Xu & Cao & 2018} and \cite[Theorem 1.1]{Xu & Cao & 2020}, a lacuna is found when the authors were trying to come up with a specific form of solution. For a transcendental entire function $f(z_1,z_2)$ with finite order satisfying the relation \beas \left(\frac{\partial f(z_1,z_2)}{\partial z_1}+if(z_1+c_1,z_2+c_2)\right)\left(\frac{\partial f(z_1,z_2)}{\partial z_1}-if(z_1+c_1,z_2+c_2)\right)=1,\eeas
the authors assumed that \beas \frac{\partial f(z_1,z_2)}{\partial z_1}+if(z_1+c_1,z_2+c_2)=e^{i p(z_1,z_2)}\eeas and \beas \frac{\partial f(z_1,z_2)}{\partial z_1}-if(z_1+c_1,z_2+c_2)=e^{-i p(z_1,z_2)},\eeas
whereas the correct assumption will be \beas \frac{\partial f(z_1,z_2)}{\partial z_1}+if(z_1+c_1,z_2+c_2)=\beta_1e^{\alpha p(z_1,z_2)}\eeas and \beas \frac{\partial f(z_1,z_2)}{\partial z_1}-if(z_1+c_1,z_2+c_2)=\beta_2e^{-\alpha p(z_1,z_2)},\eeas where $\alpha$, $\beta_1$ and $\beta_2$ are complex constants in one variable such that $\beta_1\beta_2=1$, and $p(z_1,z_2)$ is a polynomial. Therefore, in order to find out the complete characterization of the solutions of \eqref{e2.3} that Xu and Cao considered in their investigations, it is required to answer the following question.
\begin{ques}\label{q1.1}
	What could be the precise form of solutions of the Fermat-type partial differential-difference equation $(\ref{e2.3})$?
\end{ques}
\par We prove the following result which answers Question \ref{q1.1} completely.
\begin{theo}\label{t3}
Let $c=(c_1,c_2)$ be a constant in $\mathbb{C}^2$. Then any transcendental entire solution with finite order of the partial differential-difference equation $(\ref{e2.3})$ is of the form \beas f(z_1,z_2)=\frac{A}{2i}\left(\beta_1e^{\alpha\left(\displaystyle\frac{i}{\alpha}Az_1+Bz_2+h(z_2)\right)}-\beta_2e^{-\alpha\left(\displaystyle\frac{i}{\alpha}Az_1+Bz_2+h(z_2)\right)}\right),\eeas where $A$, $B$, $\alpha$, $\beta_1$, $\beta_2$ are complex constants in one variable with $A^2=1$, $\beta_1\beta_2=1$, and $h$ is a polynomial in $z_2$ satisfying the relation $h(z_2)-h(z_2+c_2)=\frac{i}{\alpha}Ac_1+Bc_2+\frac{i}{\alpha}\log A$, and $Ae^{iAc_1+\alpha Bc_2}=1$. In the special case whenever $c_2\neq0$, $h$ reduces to a constant and the form of the solution is \beas f(z_1,z_2)=\frac{A}{2i}\left(\beta_1e^{\alpha\left(\displaystyle\frac{i}{\alpha}Az_1+Bz_2+\text{constant}\right)}-\beta_2e^{-\alpha\left(\displaystyle\frac{i}{\alpha}Az_1+Bz_2+\text{constant}\right)}\right).\eeas 
\end{theo}
\begin{rem}
Let $A=1$, $\alpha=i$, $\beta_1=\beta_2=1$. Choose $c_1$, $c_2$ in $\mathbb{C}$ such that $Ae^{i(Ac_1+Bc_2)}=1$, and $h$ is a polynomial in $z_2$ only such that $h(z_2)=h(z_2+c_2)$. Then from Theorem \ref{t3}, it is clear that the solution of the partial differential-difference equation \eqref{e2.3} is \beas f(z_1,z_2)=\sin(Az_1+Bz_2+h(z_2)).\eeas 

Hence, Theorem \ref{t3} is a significant improvement of the result of Xu and Cao \cite[Theorem 1.2]{Xu & Cao & 2018}.
\end{rem}
\begin{exm}
Let $\beta_1=\beta_2=1$, $\alpha=1$, $A=B=1$ and $h(z_2)=7$, where $h$ is defined in Theorem \ref{t3}. Choose $c=(c_1,c_2)$ in $\mathbb{C}^2$ such that $e^{ic_1+c_2}=1$, $c_2\neq0$. Then it can be easily verified that \beas f(z_1,z_2)=-i\sinh(iz_1+z_2+7)\eeas is a solution of the Fermat-type partial differential-difference equation \eqref{e2.3}.
\end{exm}
\begin{exm}
Let $\beta_1=\beta_2=1$, $\alpha=1$, $A=1$, $B=3$ and $h(z_2)=z_2^3+z_2^2+13$, where $h$ is defined in Theorem \ref{t3}. Choose $c=(c_1,0)$ in $\mathbb{C}^2$ such that $e^{ic_1}=1$. Then it can be easily verified that \beas f(z_1,z_2)=-i\sinh(iz_1+3z_2+z_2^3+z_2^2+13)\eeas is a solution of the Fermat-type partial differential-difference equation \eqref{e2.3}.
\end{exm}

Next, we are interested to find solutions of Fermat-type partial differential-difference equations. Henceforth, we consider the following equations
\bea\label{e2.4}  \left(\frac{\partial f(z_1,z_2)}{\partial z_1}\right)^2+(f(z_1+c_1,z_2+c_2)-f(z_1,z_2))^2=1\eea and \bea\label{e2.5} f^2(z_1,z_2)+P^2(z_1,z_2)\left(\frac{\partial f(z_1+c_1,z_2+c_2)}{\partial z_1}-\frac{\partial f(z_1,z_2)}{\partial z_1}\right)^2=1,\eea where $P(z_1,z_2)$ is a non-zero polynomial in $\mathbb{C}^2$. For solutions of the system of Fermat-type partial differential-difference equations analogue to \eqref{e2.4}, we refer to the article \cite{Xu-Liu-Li-JMAA-2020} and references therein.\vspace{1.2mm}
\par We prove the following result finding the precise form of the solutions of \eqref{e2.4}. 
\begin{theo}\label{t4}
Let $c=(c_1,c_2)$ be a constant in $\mathbb{C}^2$. Then any transcendental entire solution with finite order of the Fermat-type partial differential-difference equation $(\ref{e2.4})$ is of the form \beas f(z_1,z_2)=-\frac{1}{4i}\left(\beta_1e^{\alpha\left(\displaystyle\frac{-2i}{\alpha}z_1+Bz_2+h(z_2)\right)}-\beta_2e^{-\alpha\left(\displaystyle\frac{-2i}{\alpha}z_1+Bz_2+h(z_2)\right)}\right),\eeas where $B$, $\alpha$, $\beta_1$, $\beta_2$ are complex constants in one variable with $\beta_1\beta_2=1$, and $h$ is a polynomial in $z_2$ only satisfying the relation 
\beas h(z_2)-h(z_2+c_2)=-\frac{2ic_1}{\alpha}+Bc_2+\frac{1}{\alpha}\log (-1).\eeas
In the special case when $c_2\neq 0$, $h$ must be constant and the solution is of the form 
\beas f(z_1,z_2)=-\frac{1}{4i}\left(\beta_1e^{\alpha\left(\displaystyle\frac{-2i}{\alpha}z_1+Bz_2+C)\right)}-\beta_2e^{-\alpha\left(\displaystyle\frac{-2i}{\alpha}z_1+Bz_2+C\right)}\right),\eeas where $B$ and $C$ are complex constant such that $e^{2ic_1-\alpha Bc_2}=-1$. 
\end{theo}

\par Next, we list some examples to exhibit the existence of solution of the Fermat-type partial differential-difference equation \eqref{e2.4}.
\begin{exm}
Let $\beta_1=\beta_2=1$, $\alpha=i$, $B=3$. Choose $c=(c_1,0)$ in $\mathbb{C}^2$ such that $e^{2ic_1}=-1$. Then from Theorem \ref{t4}, it is clear that \beas f(z_1,z_2)=-\frac{1}{2}\sin(-2z_1+3z_2+z_2^5+z_2^3+1)\eeas is a solution of the Fermat-type partial differential-difference equation \eqref{e2.4}.
\end{exm}
\begin{exm}
Let $\beta_1=\beta_2=1$, $\alpha=i$, $B=5$. Choose $c=(c_1,c_2)$ in $\mathbb{C}^2$ such that $e^{2ic_1-5c_2}=-1$, $c_2\neq0$. Then from Theorem \ref{t4}, it is clear that \beas f(z_1,z_2)=-\frac{1}{2}\sin(-2z_1+5z_2+11)\eeas is a solution of the Fermat-type partial differential-difference equation \eqref{e2.4}.
\end{exm}
\begin{exm}
Let $\beta_1=\beta_2=1$, $\alpha=1$, $B=1$ and $h(z_2)=10$, where $h$ is defined in Theorem \ref{t4}. Choose $c=(c_1,c_2)$ in $\mathbb{C}^2$ such that $e^{-(-2ic_1+c_2)}=-1$, $c_2\neq0$. Then one can easily verify that \beas f(z_1,z_2)=-\frac{1}{2i}\sinh(-2iz_1+z_2+10)\eeas is a solution of the Fermat-type partial differential-difference equation \eqref{e2.4}.
\end{exm}
\begin{exm}
Let $\beta_1=\beta_2=1$, $\alpha=1$, $B=15$ and $h(z_2)=z_2^6+z_2^4+z_2^2+1$, where $h$ is defined in Theorem \ref{t4}. Choose $c=(c_1,0)$ in $\mathbb{C}^2$ such that $e^{2ic_1}=-1$. Then one can easily verify that \beas f(z_1,z_2)=-\frac{1}{2i}\sinh(-2iz_1+15z_2+z_2^6+z_2^4+z_2^2+1)\eeas is a solution of the Fermat-type partial differential-difference equation \eqref{e2.4}.
\end{exm}

\par We prove the next result finding the precise form of the solutions of \eqref{e2.5}.
\begin{theo}\label{t5} 
Let $c=(c_1,c_2)$ be a constant in $\mathbb{C}^2$. Then any transcendental entire solution with finite order of the Fermat-type partial differential-difference equation $(\ref{e2.5})$ is of the form \beas f(z_1,z_2)=\frac{1}{2}\left(\beta_1e^{\alpha (Az_1+Bz_2+C)}+\beta_2e^{-\alpha (Az_1+Bz_2+C)}\right),\eeas where $A$, $B$, $C$ and $\alpha$ are complex constants in one variable such that $e^{\alpha(Ac_1+Bc_2)}=-1$, and $P(z_1,z_2)$ reduces to a constant $-1/2i\alpha A$.
\end{theo}

\par Next, we list some examples to exhibit the existence of solution of the Fermat-type partial differential-difference equation \eqref{e2.5}.
\begin{exm}
Let $\beta_1=\beta_2=1$, $\alpha=i$, $A=2$, $B=5$, $C=12$ and $P(z_1,z_2)=1/4$. Choose $c=(c_1,c_2)$ in $\mathbb{C}^2$ such that $e^{i(2c_1+5c_2)}=-1$. Then from Theorem \ref{t5}, it is clear that \beas f(z_1,z_2)=\cos(2z_1+5z_2+12)\eeas is a solution of the Fermat-type partial differential-difference equation \eqref{e2.5}. 
\end{exm}
\begin{exm}
Let $\beta_1=\beta_2=1$, $\alpha=1$, $A=B=C=1$ and $P(z_1,z_2)=-1/2i$. Choose $c=(c_1,c_2)$ in $\mathbb{C}^2$ such that $e^{c_1+c_2}=-1$. Then from Theorem \ref{t5}, it is clear that \beas f(z_1,z_2)=\cosh(z_1+z_2+1)\eeas is a solution of the Fermat-type partial differential-difference equation \eqref{e2.5}. 
\end{exm}

\section{\textbf{Key Lemmas}} 
In this section, we present some necessary lemmas which will play key role to prove the main results of this paper.
\begin{lem}\label{lem3.1}\cite{Hu & Li & Yang & 2003}
Let $f_j\not\equiv0$ $(j=1,2\ldots,m;\; m\geq 3)$ be meromorphic functions on $\mathbb{C}^{n}$ such that $f_1,\ldots, f_{m-1}$ are not constants, $f_1+f_2+\cdots+ f_m=1$ and such that
\beas \sum_{j=1}^{m}\left\{N_{n-1}\left(r,\frac{1}{f_j}\right)+(m-1)\ol N(r,f_j)\right\}< \lambda T(r,f_j)+O(\log^{+}T(r,f_j))\eeas
holds for $j=1,\ldots, m-1$ and all $r$ outside possibly a set with finite logarithmic measure, where $\lambda < 1$ is a positive number. Then $f_m = 1$.	
\end{lem}
\begin{lem}\label{lem3.2}\cite{Lelong & 1968, ronkin & AMS & 1971, Stoll & AMS & 1974}
For an entire function $F$ on $\mathbb{C}^n$, $F(0)\not\equiv 0$ and put $\rho(n_F)=\rho<\infty$. Then there exist a canonical function $f_F$ and a function $g_F\in\mathbb{C}^n$ such that $F(z)=f_F (z)e^{g_F(z)}$. For the special case $n=1$, $f_F$ is the canonical product of Weierstrass.
\end{lem}
\begin{lem}\label{lem3.3}\cite{P`olya & Lond & 1926}
If $g$ and $h$ are entire functions on the complex plane $\mathbb{C}$ and $g(h)$ is an entire function of finite order, then there are only two possible
cases: either
\begin{enumerate}
\item [(i)] the internal function $h$ is a polynomial and the external function $g$ is of finite order; or else
\item [(ii)] the internal function $h$ is not a polynomial but a function of finite order, and the external function $g$ is of zero order.
\end{enumerate}\end{lem}

\begin{lem}\label{lem3.4}\cite{Biancofiore & Stoll & 1981, Ye & 1995}
Let $f$ be a non-constant meromorphic function on $\mathbb{C}^n$ and let $I=(\alpha_1,\ldots,\alpha_n)$ be a multi-index with length $|I|=\sum_{j=1}^{n}\alpha_j$. Assume that
$T(r_0,f)\geq e$ for some $r_0$. Then
\beas m\left(r,\frac{\partial^{I}f}{f}\right)=S(r,f)\eeas holds for all $r\geq r_0$ outside a set $E\subset
(0,+\infty) $ of finite logarithmic measure, $\displaystyle\int_{E}\frac{dt}{t}< \infty$, where $\partial^{I}f=\displaystyle\frac{\partial^{|I|}f}{\partial z_1^{\alpha_1}\ldots\partial z_2^{\alpha_2}}$.  
\end{lem}
\begin{lem}\label{lem2.5}\cite{Cao & Korhonen & 2016, Korhonen & CMFT & 2012}
Let $f$ be a non-constant meromorphic function with finite order on $\mathbb{C}^n$ such that $f(0)\neq 0,\infty$, and let $\epsilon>0$. Then for $c\in \mathbb{C}^n$,  \beas m\left(r,\frac{f(z+c)}{f(z)}\right)+m\left(r,\frac{f(z)}{f(z+c)}\right)=S(r,f)\eeas holds for all $r\geq r_0$ outside a set $E\subset
(0,+\infty) $ of finite logarithmic measure, $\displaystyle\int_{E}\frac{dt}{t}< \infty$.\end{lem}
%\begin{lem}\label{lem3.6}\cite{Hu & Yang & 2014, Li Complex Var & 1996}Let $f$ be a non-constant meromorphic function on $\mathbb{C}^n$. Take a positive integer $m$ and take polynomials of $f$ and its partial derivatives:\beas P(f)=\sum_{p\in I}f^{p_0}(\partial^{\alpha_1}f)^{p_1}\cdots (\partial^{\alpha_l}f)^{p_l},\; (p)=(p_0,p_1,\dots,p_l)\eeas\beas Q(f)=\sum_{q\in J}f^{q_0}(\partial^{\beta_1}f)^{q_1}\cdots (\partial^{\beta_l}f)^{q_l},\; (q)=(q_0,q_1,\dots,q_l)\eeas and \beas B(f)=\displaystyle\sum_{k=0}^{m}b_kf^{k},\eeas where $I$, $J$ are finite sets of distinct elements and $a_p$, $c_q$, $b_k$ are meromorphic functions on $\mathbb{C}^n$ with $b_m\not\equiv 0$. Assume that $f$ satisfies the equation $B(f)Q(f)=P(f)$ such that $P(f)$, $Q(f)$ and $B(f)$ are differential polynomials, that is, their coefficients a satisfy $m(r, a) = S(r, f)$. If $\text{deg}(P(f))\leq m=\text{deg}(B(f))$, then \beas m(r,Q(f))=S(r,f)\eeas holds for all $r$ possibly outside of a set $E$ with finite logarithmic measure.\end{lem}
\section{Proof of the main results}
\begin{proof}[\textbf{Proof of Theorem $\ref{t1}$}]
Let $P(z)=\sum_{|I|=0}^{p}a_{\alpha_1,\ldots,\alpha_n}z_1^{\alpha_1}\ldots z_n^{\alpha_n}$ and\\ \noindent $Q(z)=\sum_{|J|=0}^{q}a_{\beta_1,\ldots,\beta_n}z_1^{\beta_1}\ldots z_n^{\beta_n}$ be two polynomials in $\mathbb{C}^n$, where $I=(\alpha_1,\ldots,\alpha_n)$, $J=(\beta_1,\ldots,\beta_n)$ be two multi-index with $|I|=\sum_{j=0}^{n}\alpha_j$ and $|I|=\sum_{j=0}^{n}\beta_j$, and $\alpha_j$, $\beta_j$ are non-negative integers.\vspace{1mm}
\par Suppose $f$ is a finite order transcendental entire solution of $(\ref{e2.1})$. We write $(\ref{e2.1})$ as \bea\label{e4.1} h_1(z)h_2(z)=Q(z),\eea where $h_1(z)=f(z)+iP(z)f(z+c)$ and $h_2(z)=f(z)-iP(z)f(z+c)$.\vspace{1mm}
\par Since $Q(0)\neq0$, hence we must have $h_1(0)\neq0$ and $h_2(0)\neq0$. Now in view of $(\ref{e4.1})$ and Lemma \ref{lem3.2}, it is easy to see that 
\bea \label{e4.2} f(z)+iP(z)f(z+c)=\chi_{_1}(z)e^{\chi_{_2}(z)};\eea
\bea \label{e4.3} f(z)-iP(z)f(z+c)=\chi_{_3}(z)e^{-\chi_{_2}(z)},\eea
where $\chi_{_1}(z)$, $\chi_{_3}(z)$ are canonical functions of $h_1$ and $h_2$, respectively and $\chi_{_2}(z)$ is entire on $\mathbb{C}^n$ such that $\chi_{_1}(z)\chi_{_3}(z)=Q(z)$.\vspace{1mm}
\par Since $f$ is a transcendental entire function of finite order, it follows from $(\ref{e4.2})$ and $(\ref{e4.3})$ that $e^{\chi_{_2}(z)}$ is of finite order, and hence by Lemma \ref{lem3.3}, it is easy to see that $\chi_{_2}(z)$ is a polynomial in $\mathbb{C}^n$.\vspace{1mm}
\par Solving $(\ref{e4.2})$ and $(\ref{e4.3})$, we obtain 
\bea\label{e4.4} f(z)=\frac{1}{2}\left(\chi_{_1}(z)e^{\chi_{_2}(z)}+\chi_{_3}(z)e^{-\chi_{_2}(z)}\right);\eea
\bea\label{e4.5} f(z+c)=\frac{1}{2iP(z)}\left(\chi_{_1}(z)e^{\chi_{_2}(z)}-\chi_{_3}(z)e^{-\chi_{_2}(z)}\right).\eea
Combining $(\ref{e4.4})$ and $(\ref{e4.5})$  yields that \beas f(z+c)&=&\frac{1}{2}\left(\chi_{_1}(z+c)e^{\chi_{_2}(z+c)}+\chi_{_3}(z+c)e^{-\chi_{_2}(z+c)}\right)\\&=&\frac{1}{2iP(z)}\left(\chi_{_1}(z)e^{\chi_{_2}(z)}-\chi_{_3}(z)e^{-\chi_{_2}(z)}\right).\eeas 
Therefore, a simple computation shows that \bea\label{e4.6} \frac{iP(z)\chi_{_1}(z+c)e^{\chi_{_2}(z+c)+\chi_{_2}(z)}}{-\chi_{_3}(z)}&+&\frac{iP(z)\chi_{_1}(z+c)e^{\chi_{_2}(z)-\chi_{_2}(z+c)}}{-\chi_{_3}(z)}\nonumber\\&&+\frac{\chi_{_1}(z)e^{2\chi_{_2}(z)}}{\chi_{_3}(z)}=1.\eea 

\par Now $\chi_{_2}(z)$ being a non-constant polynomial, it is easy to see that both 
\beas \frac{iP(z)\chi_{_1}(z+c)e^{\chi_{_2}(z+c)+\chi_{_2}(z)}}{-\chi_{_3}(z)}\;\; \mbox{and}\;\; \frac{\chi_{_1}(z)e^{2\chi_{_2}(z)}}{\chi_{_3}(z)} \eeas are non-constants.\vspace{1.2mm}

\par Since $T(r,P)=o(T(r,f))$ and $T(r,Q)=o(T(r,f))$, by Lemma \ref{lem3.1}, we obtain  \bea\label{e4.7}iP(z)\chi_{_3}(z+c)e^{\chi_{_2}(z)-\chi_{_2}(z+c)}=-\chi_{_3}(z).\eea 

\par Thus, we may write $\chi_{_2}(z)=L(z)+b$, where $L(z)=a_1z_1+a_2z_2+\cdots+a_nz_n$, $a_i,b\in \mathbb{C}$, $1\leq i\leq n$. Therefore, $(\ref{e4.7})$ reduces to the following form\bea\label{e4.8} iP(z)\chi_{_3}(z+c)e^{-\sum_{i=1}^{n}a_ic_i}=-\chi_{_3}(z).\eea 

\par Similarly, from $(\ref{e4.6})$, we easily obtain \bea\label{e4.9} iP(z)\chi_{_1}(z+c)e^{\sum_{i=1}^{n}a_ic_i}=\chi_{_1}(z).\eea 

\par Combining $\ref{e4.8}$ and $(\ref{e4.9})$ yields  \beas P^2(z)\chi_{_1}(z+c)\chi_{_3}(z+c)=\chi_{_1}(z)\chi_{_3}(z)\eeas which can be written as \beas P^2(z)Q(z+c)=Q(z).\eeas 

\par Since $Q(z)$ is a non-zero polynomial, we must have $P^2(z)=1$ and hence $P(z)=\pm1$. Consequently, $Q(z)$ reduces to a constant $q$, say. Hence, $(\ref{e2.1})$ takes the form $g^2(z)+g^2(z+c)=1$, where $g(z)=f(z)/\sqrt{q}$. Therefore, with the help of \cite[Theorem 1.4]{Xu & Cao & 2018}, we easily obtain conclusion of the theorem.
\end{proof}
\begin{proof}[\textbf{Proof of Theorem $\ref{t2}$}] We discuss the whole proof into the following two cases.\vspace{1.2mm}
\noindent\textbf{Case 1:} Let $m>n$. Since $f$ is entire, using Lemma \ref{lem3.4}, we obtain
\beas T(r, f(z_1,z_2))&=&m(r,z_1,z_2)\\&\leq& m\left(r,\frac{f(z_1,z_2)}{f(z_1+c_1,z_2+c_2)}\right)+m(r,f(z_1+c_1,z_2+c_2))+\log 2\\&\leq& m(r,f(z_1+c_1,z_2+c_2))+S(r,f)\\&\leq& T(r,f(z_1+c_1,z_2+c_2))+S(r,f).\eeas

\par By the Mohon'ko theorem for functions in $ \mathbb{C}^2 $ (see \cite{Hu & Complex Var & 1995}), an easy and straight forward computation shows that
\beas mT(r,f(z_1,z_2))&\leq& mT(r,f(z_1+c_1,z_2+c_2))+S(r,f)\\&=& T(r,f^m(z_1+c_1,z_2+c_2))+S(r,f)\\&\leq& T(r,(\partial^{I}f+\partial^{J}f)^n-1)+S(r,f)\\&\leq& nT(r,\partial^{I}f+\partial^{J}f)+S(r,f)\\&=&nm(r,\partial^{I}f+\partial^{J}f)+S(r,f)\\&\leq& nm\left(r,\frac{\partial^{I}f+\partial^{J}f}{f}\right)+nm(r,f(z_1,z_2))+S(r,f)\\&\leq& nm\left(r,\frac{\partial^{I}f}{f}\right)+m\left(r, \frac{\partial^{J}f}{f}\right)+nm(r,f(z_1,z_2))+S(r,f)\\&\leq& nT(r,f(z_1,z_2))+S(r,f).\eeas 

\par Therefore, we must have \beas (n-m)T(r,f(z_1,z_2))\leq S(r,f),\eeas which contradicts to the fact that $f$ is transcendental and also $m>n$.\vspace{1mm}
\par \textbf{Case 2:} Suppose $n>m\geq2$. In this case, it is easy to see  that $1/m+1/n<1$ and $m>n/(n-1)$. Let $d_1, d_2, \ldots, d_n$ are the roots of the equation $w^n=1$.  Then by the Second Fundamental Theorem of Nevalinna, and the equation $(\ref{e22.5})$, a simple computation shows that \beas (n-1)T(r,\partial^{I}f+\partial^{J}f)&\leq& \ol N(r,\partial^{I}f+\partial^{J}f)+\sum_{j=1}^{n}\ol N\left(r,\frac{1}{\partial^{I}f+\partial^{J}f-d_j}\right)\\&&+S(r,\partial^{I}f+\partial^{J}f)\\&\leq& \ol N\left(r,\frac{1}{(\partial^{I}f+\partial^{J}f)^n-1}\right)+S(r,\partial^{I}f+\partial^{J}f)\\&\leq& \ol N\left(r,\frac{1}{f(z_1+c_1,z_2+c_2)}\right)+S(r,\partial^{I}f+\partial^{J}f)\\&\leq& T(r,f(z_1+c_1,z_2+c_2))+S(r,f).\eeas

\par On the other hand, by Mohon'ko theorem for the functions in $ \mathbf{C}^2 $ \cite{Hu & Complex Var & 1995} and \eqref{e22.5}, we obtain
\beas mT(r,f(z_1+c_1,z_2+c_2))&=&T(r,f^m(z_1+c_1,z_2+c_2))+S(r,f)\\&=& T(r,(\partial^{I}f+\partial^{J}f)^n-1)+S(r,f)\\&=&nT(r,\partial^{I}f+\partial^{J}f)+S(r,f).\eeas 

\par Thus we have \beas mT(r,f(z_1+c_1,z_2+c_2))\leq \frac{n}{n-1}T(r,f(z_1+c_1,z_2+c_2))+S(r,f)\eeas and this can be written as \beas \left(m-\frac{n}{n-1}\right)T(r,f(z_1+c_1,z_2+c_2))\leq S(r,f),\eeas which is not possible since $f$ is transcendental and $m>{n}/{(n-1)}$. This completes the proof.
\end{proof}
\begin{proof}[\textbf{Proof of Theorem $\ref{t3}$}]
First we assume that $f$ is a finite order transcendental entire solution of $(\ref{e2.3})$. Rewriting equation $(\ref{e2.3})$, we obtain
\beas \left(\frac{\partial f(z_1,z_2)}{\partial z_1}+if(z_1+c_1,z_2+c_2)\right)\left(\frac{\partial f(z_1,z_2)}{\partial z_1}-if(z_1+c_1,z_2+c_2)\right)=1.\eeas 

\par From this equation, it is easy to see that both $\frac{\partial f(z_1,z_2)}{\partial z_1}+if(z_1+c_1,z_2+c_2)$ and $\frac{\partial f(z_1,z_2)}{\partial z_1}-if(z_1+c_1,z_2+c_2)$ do not have any zeros in $\mathbb{C}^2$.\vspace{1mm} \par Hence, in view of Lemma \ref{lem3.2}, we may assume  
\beas \frac{\partial f(z_1,z_2)}{\partial z_1}+if(z_1+c_1,z_2+c_2)=\beta_1e^{\alpha p(z_1,z_2)}\eeas and \beas \frac{\partial f(z_1,z_2)}{\partial z_1}-if(z_1+c_1,z_2+c_2)=\beta_2e^{-\alpha p(z_1,z_2)},\eeas where $p(z_1,z_2)$ is a non-constant entire function in $\mathbb{C}^2$, $\alpha$, $\beta_1$, $\beta_2$ are complex constants such that $\beta_1\beta_2=1$.

From these last two equations, a simple computation shows that
\bea\label{e4.10} \frac{\partial f(z_1,z_2)}{\partial z_1}=\frac{\beta_1e^{\alpha p(z_1,z_2)}+\beta_2e^{-\alpha p(z_1,z_2)}}{2}\eea and
\bea\label{e4.11} f(z_1+c_1,z_2+c_2)=\frac{\beta_1e^{\alpha p(z_1,z_2)}-\beta_2e^{-\alpha p(z_1,z_2)}}{2i}.\eea 

\par Using Lemma \ref{lem3.3}, it is easy to see that $p(z_1,z_2)$ must be a non-constant polynomial in $\mathbb{C}^2$. Hence, it follows from $(\ref{e4.11})$ that $p(z_1,z_2)$ is a non-constant polynomial in $\mathbb{C}^2$.\vspace{1.2mm}
\par From $(\ref{e4.10})$ and $(\ref{e4.11})$, we obtain the following \beas \frac{\partial f(z_1+c_1,z_2+c_2)}{\partial z_1}&=&\frac{\beta_1e^{\alpha p(z_1+c_1,z_2+c_2)}+\beta_2e^{-\alpha p(z_1+c_1,z_2+c_2)}}{2}\\&=& \alpha\frac{\partial p(z_1,z_2)}{\partial z_1}\frac{\beta_1e^{\alpha p(z_1,z_2)}+\beta_2e^{-\alpha p(z_1,z_2)}}{2i}.\eeas

Therefore, an elementary computation shows that 
\bea\label{e4.12}  \frac{\beta_1\alpha}{i\beta_2}\frac{\partial p(z_1,z_2)}{\partial z_1}e^{\alpha p(z_1,z_2)+ip(z_1+c_1,z_2+c_2)}&+& \frac{\alpha}{i}\frac{\partial p(z_1,z_2)}{\partial z_1}e^{\alpha p(z_1+c_1,z_2+c_2)-ip(z_1,z_2)}\nonumber\\&&-\frac{\beta_1}{\beta_2}e^{2\alpha p(z_1+c_1,z_2+c_2)}=1.\eea

From the above equation, we see that $\frac{\partial p(z_1,z_2)}{\partial z_1}$ is non-zero polynomial. For otherwise, it follows from $(\ref{e4.12})$ that $e^{2\alpha p(z_1+c_1,z_2+c_2)}=-\beta_2^2=\text{constant}$, which implies that $p(z_1,z_2)$ is constant, a contradiction. Hence, both $e^{2\alpha p(z_1+c_1,z_2+c_2)}$ and $ \frac{\beta_1\alpha}{i\beta_2}\frac{\partial p(z_1,z_2)}{\partial z_1}e^{\alpha p(z_1,z_2)+\alpha p(z_1+c_1,z_2+c_2)}$ are non-constants.

Also notice that 
\beas N(r,e^{\alpha 2p(z_1+c_1,z_2+c_2)})=N\left(r,\frac{1}{e^{\alpha 2p(z_1+c_1,z_2+c_2)}}\right)=S(r,f),\eeas

\beas&& N\left(r,\frac{\beta_1\alpha}{i\beta_2}\frac{\partial p(z_1,z_2)}{\partial z_1}e^{\alpha p(z_1,z_2)+\alpha p(z_1+c_1,z_2+c_2)}\right)\\&=&N\left(r,\frac{i\beta_2}{\frac{\beta_1\alpha \partial p(z)}{\partial z_1}e^{\alpha p(z_1,z_2)+\alpha p(z_1+c_1,z_2+c_2)}}\right)=S(r,f).\eeas

\par Therefore, by Lemma \ref{lem3.1}, we obtain 
\beas -i\alpha\frac{\partial p(z_1,z_2)}{\partial z_1}e^{\alpha p(z_1+c_1,z_2+c_2)-\alpha p(z_1,z_2)}=1,\eeas which implies that \bea\label{e4.13} -i\alpha\frac{\partial p(z_1,z_2)}{\partial z_1}=e^{\alpha p(z_1,z_2)-\alpha p(z_1+c_1,z_2+c_2)}.\eea 

\par From $(\ref{e4.12})$ and $(\ref{e4.13})$, it is easy to obtain  \bea\label{e4.15} -i\alpha\frac{\partial p(z_1,z_2)}{\partial z_1}=e^{\alpha p(z_1+c_1,z_2+c_2)-\alpha p(z_1,z_2)}.\eea
\par We observe that L.H.S. and R.H.S. of \eqref{e4.13} are respectively, polynomial (non-transcendental)  and transcendental entire function. Therefore, the only possibility is that $\alpha p(z_1,z_2)-\alpha p(z_1+c_1,z_2+c_2)$ and hence $e^{\alpha p(z_1,z_2)-\alpha p(z_1+c_1,z_2+c_2)}$ is a constant.\vspace{1.2mm}
Assume that \bea\label{e4.15a} -i\alpha\frac{\partial p(z_1,z_2)}{\partial z_1}=e^{\alpha p(z_1,z_2)-\alpha p(z_1+c_1,z_2+c_2)}=A,\eea where $A$ is a non-zero complex constant. Then from \eqref{e4.15} , we obtain
\bea\label{e4.15b} -i\alpha\frac{\partial p(z_1,z_2)}{\partial z_1}=e^{\alpha p(z_1+c_1,z_2+c_2)-\alpha p(z_1,z_2)}=\frac{1}{A}.\eea

\par From \eqref{e4.15a} and \eqref{e4.15b}, it can be easily seen that $A^2=1$. Again in view of \eqref{e4.15a}, we may assume that $p(z_1,z_2)=iAz_1/\alpha+g(z_2)$, where $g(z_2)$ is a polynomial in $z_2$ only.\vspace{1mm}

From \eqref{e4.15a}, we get \beas p(z_1,z_2)-p(z_1+c_1,z_2+c_2)=\frac{1}{\alpha}\log A.\eeas 

\par This implies that \beas g(z_2)-g(z_2+c_2)=\frac{i}{\alpha}Ac_1+\frac{1}{\alpha}\log A.\eeas
 
\par We may write $g(z_2)=Bz_2+h(z_2)$, where $B$ is a complex constant in one variable and $h$ is a polynomial in one variable $z_2$ of degree greater than one. Then we easily obtain
\bea\label{e4.15c} h(z_2)-h(z_2+c_2)=\frac{i}{\alpha}Ac_1+Bc_2+\frac{1}{\alpha}\log A.\eea

Now if $c_2=0$, then from \eqref{e4.15c}, we obtain $Ae^{iAc_1}=1$. If $c_2\neq0$, then from \eqref{e4.15c}, it is clear that $h(z_2)$ must be constant, and in that case we obtain $Ae^{iAc_1+\alpha Bc_2}=1$. Since $A^2=1$, we have $iAc_1+\alpha Bc_2=k\pi$, where $k$ is an integer.\vspace{1mm}

\par Hence, keeping in view of \eqref{e4.11}, \eqref{e4.15c}, we get after simple calculation
\beas
f(z_1,z_2)=\frac{A}{2i}\left(\beta_1e^{\alpha\left(\displaystyle\frac{i}{\alpha}Az_1+Bz_2+h(z_2)\right)}-\beta_2e^{-\alpha\left(\displaystyle\frac{i}{\alpha}Az_1+Bz_2+h(z_2)\right)}\right).\eeas

In the special case whenever $c_2\neq0$, $h$ reduces to a constant and the form of the solution is 
\beas
f(z_1,z_2)=\frac{A}{2i}\left(\beta_1e^{\alpha\left(\displaystyle\frac{i}{\alpha}Az_1+Bz_2+\text{constant}\right)}-\beta_2e^{-\alpha\left(\displaystyle\frac{i}{\alpha}Az_1+Bz_2+\text{constant}\right)}\right).\eeas

\par This completes the proof.
\end{proof}
\begin{proof}[\textbf{Proof of Theorem $\ref{t4}$}]
Assume that $f$ is a finite order transcendental entire solution of $(\ref{e2.4})$. Let us denote $f(z_1+c_1,z_2+c_2)-f(z_1,z_2)$ by $\Delta_cf(z_1,z_2)$. Rewriting equation $(\ref{e2.4})$ we have
\beas \left(\frac{\partial f(z_1,z_2)}{\partial z_1}+i\Delta_cf(z_1,z_2)\right)\left(\frac{\partial f(z_1,z_2)}{\partial z_1}-i\Delta_cf(z_1,z_2)\right)=1.\eeas 

\par Then by the same argument in the proof of Theorem \ref{t3}, we see that 
\bea\label{e4.16} \frac{\partial f(z_1,z_2)}{\partial z_1}=\frac{\beta_1e^{\alpha p(z_1,z_2)}+\beta_2e^{-\alpha p(z_1,z_2)}}{2}\eea and
\bea\label{e4.17} f(z_1+c_1,z_2+c_2)-f(z_1,z_2)=\frac{\beta_1e^{\alpha p(z_1,z_2)}-\beta_2e^{-\alpha p(z_1,z_2)}}{2i}.\eea

\par Differentiating $(\ref{e4.17})$ partially with respect to $z_1$, we obtain
\bea\label{e4.18} \frac{\partial f(z_1+c_1,z_2+c_2)}{\partial z_1}=\left(1-i\alpha\frac{\partial p(z_1,z_2)}{\partial z_1}\right)\frac{\beta_1e^{\alpha p(z_1,z_2)}+\beta_2e^{-\alpha p(z_1,z_2)}}{2}.\eea

\par Now keeping in view of \eqref{e4.16} and \eqref{e4.18}, we obtain after simple calculation that \bea\label{e4.19} &&\frac{\beta_1}{\beta_2}\left(1-i\alpha\frac{\partial p(z_1,z_2)}{\partial z_1}\right)e^{\alpha p(z_1,z_2)+\alpha p(z_1+c_1,z_2+c_2)}\nonumber\\&&+\left(1-i\alpha\frac{\partial p(z_1,z_2)}{\partial z_1}\right)e^{-\alpha p(z_1,z_2)+\alpha p(z_1+c_1,z_2+c_2)}\nonumber\\&&-\frac{\beta_1}{\beta_2}e^{2\alpha p(z_1+c_1,z_2+c_2)}=1.\eea

From the above equation, we see that $\left(1-i\alpha\frac{\partial p(z_1,z_2)}{\partial z_1}\right)$ is non-zero polynomial. For otherwise, it follows from \eqref{e4.19} that $e^{2\alpha p(z_1+c_1,z_2+c_2)}=-\beta_1^2=\text{constant}$, which implies that $p(z_1,z_2)$ is constant, a contradiction. Hence, both $-\frac{\beta_1}{\beta_2}e^{2\alpha p(z_1+c_1,z_2+c_2)}$ and $ \frac{\beta_1}{\beta_2}\left(1-i\alpha\frac{\partial p(z_1,z_2)}{\partial z_1}\right)\frac{\partial p(z_1,z_2)}{\partial z_1}e^{\alpha p(z_1,z_2)+\alpha p(z_1+c_1,z_2+c_2)}$ are non-constants.

 We also observe that 
\beas N\left(r,e^{\alpha2p(z_1+c_1,z_2+c_2)}\right)=N\left(r,\frac{1}{e^{\alpha2p(z_1+c_1,z_2+c_2)}}\right)=S(r,f),\eeas

\beas&& N\left(r,\left(1-i\alpha\frac{\partial p(z_1,z_2)}{\partial z_1}\right)e^{\alpha ip(z_1,z_2)+\alpha p(z_1+c_1,z_2+c_2)}\right)\\&&=N\left(r,\frac{1}{\left(1-i\alpha\frac{\partial p(z_1,z_2)}{\partial z_1}\right)e^{\alpha ip(z_1,z_2)+\alpha p(z_1+c_1,z_2+c_2)}}\right)=S(r,f)\eeas and

\beas&& N\left(r,\left(1-i\alpha\frac{\partial p(z_1,z_2)}{\partial z_1}\right)e^{-\alpha ip(z_1,z_2)+\alpha p(z_1+c_1,z_2+c_2)}\right)\\&&=N\left(r,\frac{1}{\left(1-i\alpha\frac{\partial p(z_1,z_2)}{\partial z_1}\right)e^{-\alpha ip(z_1,z_2)+\alpha p(z_1+c_1,z_2+c_2)}}\right)=S(r,f).\eeas

By Lemma \ref{lem3.1}, we easily obtain 
\bea\label{e4.20} \left(1-i\alpha\frac{\partial p(z_1,z_2)}{\partial z_1}\right)e^{ip(z_1+c_1,z_2+c_2)-ip(z_1,z_2)}=1,\eea
which implies that \bea\label{e4.21} 1-i\alpha\frac{\partial p(z_1,z_2)}{\partial z_1}=e^{ip(z_1,z_2)-ip(z_1+c_1,z_2+c_2)}.\eea

\par Again, using \eqref{e4.19} and \eqref{e4.20}, we obtain \bea\label{e4.22} 1-i\alpha\frac{\partial p(z_1,z_2)}{\partial z_1}=e^{\alpha p(z_1+c_1,z_2+c_2)-ip(z_1,z_2)}.\eea

Now we observe that L.H.S. of $(\ref{e4.21})$ is a polynomial, whereas R.H.S. of it is transcendental. Therefore, we conclude that $\alpha p(z_1,z_2)-\alpha p(z_1+c_1,z_2+c_2)$, and thus $e^{\alpha p(z_1,z_2)-\alpha p(z_1+c_1,z_2+c_2)}$ must be constant.\vspace{1.2mm}

Let \bea\label{e4.22a}  1-i\alpha\frac{\partial p(z_1,z_2)}{\partial z_1}=e^{ip(z_1,z_2)-ip(z_1+c_1,z_2+c_2)}=A,\eea where $A$ is a non-zero complex constant in one variable.\vspace{1.2mm}

Then from \eqref{e4.22}, we get 
\bea\label{e4.22b} 1-i\alpha\frac{\partial p(z_1,z_2)}{\partial z_1}=e^{\alpha p(z_1+c_1,z_2+c_2)-ip(z_1,z_2)}=\frac{1}{A}.\eea

Therefore, it follows from \eqref{e4.22a} and \eqref{e4.22b} that $A^2=1$.\vspace{1.2mm}

\par Again, since $1-i\alpha\frac{\partial p(z_1,z_2)}{\partial z_1}=A$, we may assume that \beas p(z_1,z_2)=\frac{A-1}{\alpha}iz_1+g(z_2),\eeas where $g$ is a polynomial in one variable $z_2$ only.\vspace{1.2mm}

\par From \eqref{e4.22a}, it follows that $p(z_1,z_2)-p(z_1+c_1,z_2+c_2)=\frac{1}{\alpha}\log A$. This implies that 
\beas g(z_2)-g(z_2+c_2)=\frac{(A-1)i}{\alpha}c_1+\frac{1}{\alpha}\log A.\eeas

Let $g(z_2)=Bz_2+h(z_2)$, where $h$ is a polynomial in $z_2$ only with degree greater than 1 or constant and $B$ is a complex number in one variable.\vspace{1.2mm}

\par Then \beas p(z_1,z_2)=\frac{(A-1)i}{\alpha}z_1+Bz_2+h(z_2)\;\; \text{and}\eeas
\bea\label{e4.22c} h(z_2)-h(z_2+c_2)=\frac{(A-1)i}{\alpha}z_1+Bz_2+\frac{1}{\alpha}\log A.\eea

If $c_2=0$, then from \eqref{e4.22c}, it follows that $Ae^{i(A-1)}=1$. If $c\neq0$, then by \eqref{e4.22c}, we conclude that $h$ must be constant and $Ae^{i(A-1)c_1+\alpha Bc_2}=1$.\vspace{1.2mm}

\par From \eqref{e4.16}, we obtain 
\beas f(z_1,z_2)=\frac{1}{2i(A-1)}\left(\beta_1e^{i(A-1)z_1+\alpha Bz_2+\alpha h(z_2)}-\beta_2e^{-(i(A-1)z_1+\alpha Bz_2+\alpha h(z_2))}\right),\eeas where $A\neq1$.\vspace{1.2mm} 

\par Since $A^2=1$, it must be $A=-1$. Therefore, the form of the function is 
\beas f(z_1,z_2)=-\frac{1}{4i}\left(\beta_1e^{\alpha\left(\displaystyle\frac{-2i}{\alpha}z_1+Bz_2+h(z_2)\right)}-\beta_2e^{-\alpha\left(\displaystyle\frac{-2i}{\alpha}z_1+Bz_2+h(z_2)\right)}\right),\eeas where $h$ satisfies the relation
\beas h(z_2)-h(z_2+c_2)=-\frac{2i}{\alpha}z_1+Bz_2+\frac{1}{\alpha}\log(-1).\eeas

In the special case whenever $c_2\neq0$, $h$ is constant by \eqref{e4.22c} and the form of the function is of the form 
\beas f(z_1,z_2)=-\frac{1}{4i}\left(\beta_1e^{\alpha\left(\displaystyle\frac{-2i}{\alpha}z_1+Bz_2+C\right)}-\beta_2e^{-\alpha\left(\displaystyle\frac{-2i}{\alpha}z_1+Bz_2+C\right)}\right),\eeas where $B$, $C$ are complex constants such that $e^{2iz_1-Bz_2}=-1$.

This completes the proof of the theorem.
\end{proof}

\begin{proof}[\textbf{Proof of Theorem $\ref{t5}$}]	First we assume that $f$ is a finite order transcendental entire solution of $(\ref{e2.5})$. Let us rewrite $(\ref{e2.5})$ as \bea \label{e3.27}h_1(z_1,z_2)h_2(z_1,z_2)=1,\eea where 
\[
\begin{cases}
	 \displaystyle h_1(z_1,z_2)=f(z_1,z_2)+iP(z_1,z_2)\left(\frac{\partial f(z_1+c_1,z_2+c_2)}{\partial z_1}-\frac{\partial f(z_1,z_2)}{\partial z_1}\right)\vspace{1.5mm}\\ \displaystyle h_2(z_1,z_2)=f(z_1,z_2)-iP(z_1,z_2)\left(\frac{\partial f(z_1+c_1,z_2+c_2)}{\partial z_1}-\frac{\partial f(z_1,z_2)}{\partial z_1}\right).
\end{cases}	
\]	

In view of $(\ref{e3.27})$ and Lemma \ref{lem3.2}, we assume that
\bea\label{e3.28} f(z_1,z_2)+iP(z_1,z_2)\left(\frac{\partial f(z_1+c_1,z_2+c_2)}{\partial z_1}-\frac{\partial f(z_1,z_2)}{\partial z_1}\right)=\beta_1e^{\alpha
H(z)}\eea and \bea\label{e3.29} f(z_1,z_2)-iP(z_1,z_2)\left(\frac{\partial f(z_1+c_1,z_2+c_2)}{\partial z_1}-\frac{\partial f(z_1,z_2)}{\partial z_1}\right)=\beta_2e^{-\alpha H(z)},\eea	where $H(z)$ is an entire functions in $\mathbb{C}^2$ and $\alpha(\neq0)$, $\beta_1$, $\beta_2$ are complex constants such that $\beta_1\beta_2=1$.\vspace{1mm}
\par Solving $(\ref{e3.28})$ and $(\ref{e3.29})$, we obtain 
\bea\label{e3.30} f(z_1,z_2)=\frac{\beta_1e^{\alpha H(z_1,z_2)}+\beta_2e^{-\alpha H(z_1,z_2)}}{2}\eea and \bea\label{e3.31} P(z_1,z_2)\left(\frac{\partial f(z_1+c_1,z_2+c_2)}{\partial z_1}-\frac{\partial f(z_1,z_2)}{\partial z_1}\right)=\frac{\beta_1e^{\alpha H(z_1,z_2)}-\beta_2e^{-\alpha H(z_1,z_2)}}{2i}.\eea

Since $f(z_1,z_2)$ is a transcendental entire function of finite order, in view of $(\ref{e3.30})$, it is easy to see that $H(z_1,z_2)$ is a polynomial in $\mathbb{C}^2$.\vspace{1mm}
\par Differentiating $(\ref{e3.30})$ partially with respect to $z_1$, we obtain
\bea\label{e3.32} \frac{\partial f(z_1,z_2)}{\partial z_1}=\frac{\alpha}{2}\frac{\partial H(z_1,z_2)}{\partial z_1}\left(\beta_1e^{\alpha H(z_1,z_2)}-\beta_2e^{-\alpha H(z_1,z_2)}\right).\eea
\par From $(\ref{e3.31})$ and $(\ref{e3.32})$, we obtain
\bea\label{e3.33}&& \frac{\left(\alpha\frac{\partial H(z_1,z_2)}{\partial z_1}+\frac{1}{iP(z_1,z_2)}\right)\beta_1}{-\alpha \beta_2\frac{\partial H(z_1+c_1,z_2+c_2)}{\partial z_1}}e^{\alpha H(z_1,z_2)+\alpha H(z_1+c_1,z_2+c_2)}\nonumber\\&&+\frac{\left(\alpha\frac{\partial H(z_1,z_2)}{\partial z_1}+\frac{1}{iP(z_1,z_2)}\right)}{\alpha \frac{\partial H(z_1+c_1,z_2+c_2)}{\partial z_1}}e^{-\alpha H(z_1,z_2)+\alpha H(z_1+c_1,z_2+c_2)}\nonumber\\&&+\frac{\beta_1}{\beta_2}e^{2\alpha H(z_1+c_1,z_2+c_2)}=1.\eea

Since $H(z_1,z_2)$ is a non-constant polynomial, it follows from \eqref{e3.33} that $\alpha\frac{\partial H(z_1,z_2)}{\partial z_1}+\frac{1}{iP(z_1,z_2)}$ can not be zero. For otherwise, from \eqref{e3.33}, we obtain $e^{2\alpha H(z_1+c_1,z_2+c_2)}=-\beta_1^2=$ constant, which implies that $H$ is constant, which is a contradiction.\vspace{1.2mm}

Therefore, both
 \beas\frac{\left(\alpha\frac{\partial H(z_1,z_2)}{\partial z_1}+\frac{1}{iP(z_1,z_2)}\right)\beta_1}{-\alpha \beta_2\frac{\partial H(z_1+c_1,z_2+c_2)}{\partial z_1}}e^{\alpha H(z_1,z_2)+\alpha H(z_1+c_1,z_2+c_2)} \;\text{and}\;\frac{\beta_1}{\beta_2}e^{\alpha2H(z_1+_1,z_2+c_2)}\eeas are not constants.\vspace{1.2mm} 
 
 \par Also in view of Section 3.11 (page no. $301$) of \cite{Hu & Li & Yang & 2003}, a simple computation shows that 
\beas &&N\left(r, \frac{\left(\alpha\frac{\partial H(z_1,z_2)}{\partial z_1}+\frac{1}{iP(z_1,z_2)}\right)\beta_1}{-\alpha \beta_2\frac{\partial H(z_1+c_1,z_2+c_2)}{\partial z_1}}e^{\alpha H(z_1,z_2)+\alpha H(z_1+c_1,z_2+c_2)}\right)\\&&=N\left(r, \frac{-\alpha \beta_2\frac{\partial H(z_1+c_1,z_2+c_2)}{\partial z_1}}{\left(\alpha\frac{\partial H(z_1,z_2)}{\partial z_1}+\frac{1}{iP(z_1,z_2)}\right)\beta_1}e^{-\alpha H(z_1,z_2)-\alpha H(z_1+c_1,z_2+c_2)}\right)=S(r,f),\eeas

\beas && N\left(r, \frac{\left(\alpha\frac{\partial H(z_1,z_2)}{\partial z_1}+\frac{1}{iP(z_1,z_2)}\right)}{\alpha \frac{\partial H(z_1+c_1,z_2+c_2)}{\partial z_1}}e^{-\alpha H(z_1,z_2)+\alpha H(z_1+c_1,z_2+c_2)}\right)\\&&=N\left(r, \frac{\alpha \frac{\partial H(z_1+c_1,z_2+c_2)}{\partial z_1}}{\left(\alpha\frac{\partial H(z_1,z_2)}{\partial z_1}+\frac{1}{iP(z_1,z_2)}\right)}e^{\alpha H(z_1,z_2)-\alpha H(z_1+c_1,z_2+c_2)}\right)=S(r,f)\eeas and 
\beas&& N\left(r, e^{2\alpha H(z_1+c_1,z_2+c_2)}\right)=N\left(r, e^{-2\alpha H(z_1+c_1,z_2+c_2)}\right)=S(r,f).\eeas

Hence, applying Lemma \ref{lem3.1}, we obtain
\bea\label{e3.34} \frac{\left(\alpha\frac{\partial H(z_1,z_2)}{\partial z_1}+\frac{1}{iP(z_1,z_2)}\right)}{\alpha \frac{\partial H(z_1+c_1,z_2+c_2)}{\partial z_1}}e^{-\alpha H(z_1,z_2)+\alpha H(z_1+c_1,z_2+c_2)}=1,\eea
which implies that 
\bea\label{e3.34a} \frac{\left(\alpha\frac{\partial H(z_1,z_2)}{\partial z_1}+\frac{1}{iP(z_1,z_2)}\right)}{\alpha \frac{\partial H(z_1+c_1,z_2+c_2)}{\partial z_1}}=e^{\alpha H(z_1,z_2)-\alpha H(z_1+c_1,z_2+c_2)}.\eea

Next, using $(\ref{e3.33})$ and $(\ref{e3.34})$, we obtain  
\bea\label{e3.35} \frac{\left(\alpha\frac{\partial H(z_1,z_2)}{\partial z_1}+\frac{1}{iP(z_1,z_2)}\right)}{\alpha \frac{\partial H(z_1+c_1,z_2+c_2)}{\partial z_1}}=e^{\alpha H(z_1+c_1,z_2+c_2)-\alpha H(z_1,z_2)}.\eea 

\par Observe that L.H.S. of \eqref{e3.34a} is rational, whereas R.H.S. of the same is transcendental entire. Therefore, $\alpha H(z_1,z_2)-\alpha H(z_1+c_1,z_2+c_2)$ and hence\\ $e^{\alpha H(z_1,z_2)-\alpha H(z_1+c_1,z_2+c_2)}$ must be constant.\vspace{1.2mm}

\par Hence, we may assume that  \bea\label{e3.36}H(z)=Az_1+Bz_2+C,\eea where $A, B, C$ are complex constants.\vspace{1.2mm}

\par Using \eqref{e3.34}, \eqref{e3.35} and \eqref{e3.36}, we obtain 
\beas (\alpha AP(z_1,z_2))^2=(\alpha AP(z_1,z_2)-i)^2,\eeas which implies that $P(z_1,z_2)={-1}/{2i\alpha A}$=constant.\vspace{1mm}

Also, from \eqref{e3.34} and \eqref{e3.36}, we obtain that $e^{\alpha(Ac_1+Bc_2)}=-1.$ This implies that $\alpha(Ac_1+Bc_2)=(2k+1)\pi i$, where $k$ is an integer.\vspace{1.2mm} 

\par Hence from $(\ref{e3.30})$, it is easy to see that \beas f(z_1,z_2)=\frac{1}{2}\left(\beta_1e^{\alpha (Az_1+Bz_2+C)}+\beta_2e^{-\alpha (Az_1+Bz_2+C)}\right).\eeas

This completes the proof of the theorem.
\end{proof}
\vspace{1.5mm}
\noindent{\bf Acknowledgment:} The authors would like to thank the referee(s) for the helpful suggestions and comments to improve the exposition of the paper.

\end{document}